\documentclass{amsart}%
\usepackage{graphicx}
\usepackage{amscd}
\usepackage{amsmath}
\usepackage{amsfonts}
\usepackage{amssymb}%
\newtheorem{theorem}{Theorem}
\theoremstyle{plain}

\newtheorem{corollary}{Corollary}

\newtheorem{definition}{Definition}
\newtheorem{example}{Example}

\newtheorem{proposition}{Proposition}

\numberwithin{equation}{section}

\begin{document}
\title{An optimal matching problem}
\author{Ivar Ekeland}
\address{University of British Columbia, Vancouver BC, \\
V6T\ 1Z2 \ Canada}
\email{ekeland@math.ubc.ca}
\thanks{This research has been supported by NSF\ grant ***. The author thanks Jim
Heckman for introducing him to the economics of hedonic pricing}
\date{August 15, 2003}
\subjclass{Primary 05C38, 15A15; Secondary 05A15, 15A18}
\keywords{optimal transportation, measure-preserving maps}

\begin{abstract}
Given two measured spaces $(X,\mu)$ and $(Y,\nu)$, and a third space $Z$,
given two functions $u(x,z)$ and $v(x,z)$, we study the problem of finding two
maps $s:X\rightarrow Z$ and $t:Y\rightarrow Z$ such that the images $s(\mu)$
and $t(\nu)$ coincide, and the integral $\int_{X} u(x,s(x))d\mu+\int_{Y}
v(y,t(y))d\nu$ is maximal. We give condition on $u$ and $v$ for which there is
a unique solution.

\end{abstract}
\maketitle

\section{The main result.}

Suppose we are given three goods, $X,Y$, and $Z$. They are not homogeneous,
but come in different qualities, $x\in X$ , $y\in Y$ and $z\in Z$. Goods $X$
and $Y$ are used for the sole purpose of producing good $Z$, which we are
interested in. To obtain one piece of good $Z$, one has to assemble one piece
of good $X$ and one piece of good $Y$. More precisely, one can obtain a piece
of quality $z$ by assembling one piece of quality $x$ and one piece of quality
$y$, yielding a benefit of $u\left(  x,z\right)  +w\left(  y,z\right)  $.
Given the distributions $\mu$ and $\nu$ of goods $X$ and $Y$, one wishes to
minimize the total benefit of production.

This translates into the following optimization problem: find maps
$s:X\rightarrow Z$ and $t:Y\rightarrow Z$ such that $s\left(  \mu\right)
=t\left(  \nu\right)  $ (this is the matching condition) and the integral
\begin{equation}
\int_{X}u\left(  x,s\left(  x\right)  \right)  d\mu+\int_{Y}w(y,t(y))d\nu
\label{e1}%
\end{equation}
is maximized.

The origins of that problem lie in the economic theory of hedonic pricing (see
\cite{Heckman} for an overview). The economic aspects will be developed in
another paper \cite{IE1}.

Mathematically speaking, this is related to the classical optimal
transportation problem (see the monographs \cite{RaRu}and
\cite{Villani} for accounts of the theory). Recall that this
problem consists in minimizing the integral
\begin{equation}
\int_{X}u\left(  x,s\left(  x\right)  \right)  d\mu\label{e2}%
\end{equation}
among all maps $s:X\rightarrow Y$ such that $s\left(  \mu\right)  =\nu$. Here
the measured spaces $\left(  X,\mu\right)  $ and $\left(  Y,\nu\right)  $ are
given, as well as the function $u:X\times Y\rightarrow R$. A seminal result by
Brenier \cite{Brenier} states that, if $X$ and $Y$ are bounded open subsets of
$R^{n}$ endowed with the Lebesgue measure, with $X$ connected, and $u\left(
x,y\right)  =\left\|  x-y\right\|  $, then there is a unique solution $s$ to
the optimal transportation problem, and $s$ is almost everywhere equal to the
gradient of a convex function.

Kantorovitch \cite{Kant} introduced into the optimal
transportation problem a duality method which will be crucial to
our proof. Instead of proving directly existence and uniqueness in
the optimal matching problem, we solve in section 3 another
optimization problem, and we will show in section 4 that it yields
the solution to the original one. This correspondence relies
heavily on an extension of the classical duality results in convex
analysis (see \cite{IETemam}). This extension has been can be
found in \cite{RaRu} and \cite{GRC}; for the reader's convenience,
we will give the main results in section 2. Finally, in section 5,
we will give some consequences of the main result.

From now on, $X\subset R^{n_{1}},Y\subset R^{n_{2}},$ and $Z\subset R^{n_{3}}
$ will be compact subsets. We are given measures $\mu$ on $X$ and $\nu$ on
$Y$, which are absolutely continuous with respect to the Lebesgue measure, and
satisfy:
\[
\mu\left(  X\right)  =\nu\left(  Y\right)  <\infty
\]

We are also given functions $u:\Omega_{1}\rightarrow R$ and $v:\Omega
_{2}\rightarrow R$, where $\Omega_{1}$ is a neighbourhood of $X\times Z$ and
$\Omega_{2}$ is a neigbourhood of $Y\times Z$. It is assumed that $u$ and $v$
are continous with respect to both variables, differentiable with respect to
$x$ and $y,$ and that the partial derivatives $D_{x}u$ and $D_{y}v $ are
continous with respect to both variables, and injective with respect to $z:$%
\begin{align}
\forall x  & \in X,\;D_{x}u\left(  x,z_{1}\right)  =D_{x}u\left(
x,z_{2}\right)  \Longrightarrow z_{1}=z_{2}\label{eq11}\\
\forall y  & \in Y,\;D_{y}v\left(  y,z_{1}\right)  =D_{y}v\left(
y,z_{2}\right)  \Longrightarrow z_{1}=z_{2}\label{eq12}%
\end{align}

The latter condition is a generalization of the classical Spence-Mirrlees
condition in the economics of assymmetric information (see \cite{Carlier}). It
is satisfied for $u\left(  x,z\right)  =\left\|  x-z\right\|  ^{\alpha}$,
provided $\alpha\neq0$ and $\alpha\neq1$.

\begin{theorem}
\label{thm1}Under the above assumptions, there exists a pair of Borelian maps
$\left(  \bar{s},\bar{t}\right)  $ with $\bar{s}\left(  \mu\right)  =\bar
{t}\left(  \nu\right)  $, such that for every $\left(  s,t\right)  $
satisfying $s\left(  \mu\right)  =t\left(  \nu\right)  $, we have:
\[
\int_{X}u\left(  x,s\left(  x\right)  \right)  d\mu-\int_{Y}v(y,t(y))d\nu
\leq\int_{X}u\left(  x,\bar{s}\left(  x\right)  \right)  d\mu-\int_{Y}%
v(y,\bar{t}(y))d\nu<\infty
\]
This solution is unique, up to equality almost eveywhere, and it is described
as follows:\ there is some Lipschitz continuous function $\bar{p}:Z\rightarrow
R$ and some negligible subsets $X_{0}\subset X$ and $Y_{0}\subset Y$ such
that, \ for every $x\notin X_{0}$ and every $y\notin Y_{0}$:
\begin{align}
\forall z  & \neq\bar{s}\left(  x\right)  ,\;\;u\left(  x,z\right)  -\bar
{p}\left(  z\right)  <u\left(  x,\bar{s}\left(  x\right)  \right)  -\bar
{p}\left(  \bar{s}\left(  x\right)  \right) \label{eq21}\\
\forall z  & \neq\bar{t}\left(  y\right)  ,\;\;v\left(  y,z\right)  -\bar
{p}\left(  z\right)  >v\left(  x,\bar{t}\left(  y\right)  \right)  -\bar
{p}\left(  \bar{t}\left(  y\right)  \right) \label{eq22}\\
& \bar{p}\;\text{\textrm{is differentiable at\ }}\bar{s}\left(  x\right)
\mathrm{\;and\;}\bar{t}\left(  y\right) \label{eq23}%
\end{align}

\end{theorem}

If in addition $u$ and $v$ are differentiable with respect to $z$, we get,
from the minimization (\ref{eq21}) and the maximization (\ref{eq22}):
\begin{align*}
D_{z}u\left(  x,\bar{s}\left(  x\right)  \right)   & =D_{z}\bar{p}\left(
\bar{s}\left(  x\right)  \right) \\
D_{z}v\left(  y,\bar{t}\left(  y\right)  \right)   & =D_{z}\bar{p}\left(
\bar{t}\left(  y\right)  \right)
\end{align*}

Set $\bar{s}\left(  \mu\right)  =\bar{t}\left(  \nu\right)  =\lambda$. It
follows from the above that, for $\lambda$-almost every $z\in Z$, there is
some $x\notin X_{0}$ and $y\notin Y_{0}$ such that $z=\bar{s}\left(  x\right)
=\bar{t}\left(  y\right)  $, and for every such $\left(  x,y\right)  $ we
have:
\[
D_{z}u\left(  x,z\right)  =D_{z}\bar{p}\left(  z\right)  =D_{z}v\left(
y,z\right)
\]

Note that there is no reason why $\lambda$ should be absolutely continuous
with respect to the Lebesgue measure.

The proof ot theorem \ref{thm1}is deferred to section 4. Meanwhile, let us
notice that we have slightly changed the formulation of the optimal matching
problem:\ by setting $w=-v$ we recover the original one. This change will
simplify future notations.

\section{Fundamentals of $u$-convex analysis.}

In this section, we basically follow Carlier \cite{Carlier}.

\subsection{$u$-convex functions.}

We will be dealing with function taking values in $\mathbb{R\cup}\left\{
+\infty\right\}  $.

A function $f:X\rightarrow\mathbb{R\cup}\left\{  +\infty\right\}  $ will be
called $u$\emph{-convex} iff there exists a non-empty subset $A\subset
Z\times\mathbb{R}$ such that:
\begin{equation}
\forall x\in X,\;\;f\left(  x\right)  =\sup_{\left(  z,\alpha\right)  \in
A}\left\{  u\left(  x,z\right)  +a\right\} \label{u-con}%
\end{equation}

A function $p:Z\rightarrow\mathbb{R\cup}\left\{  +\infty\right\}  $ will be
called $u$\emph{-convex} iff there exists a non-empty subset $B\subset
X\times\mathbb{R}$ such that:
\begin{equation}
p\left(  z\right)  =\sup_{\left(  x,b\right)  \in B}\left\{  u\left(
x,z\right)  +b\right\} \label{u-conv*}%
\end{equation}

\subsection{Subconjugates}

Let $f:X\rightarrow\mathbb{R\cup}\left\{  +\infty\right\}  $, not identically
$\left\{  +\infty\right\}  $, be given. We define its \emph{subconjugate}
$f^{\sharp}:Z\rightarrow\mathbb{R\cup}\left\{  +\infty\right\}  $ by:
\begin{equation}
f^{\sharp}\left(  z\right)  =\sup_{x}\left\{  u\left(  x,z\right)  -f\left(
x\right)  \right\} \label{Fen}%
\end{equation}

It follows from the definitions that $f^{\sharp}$ is a $u$-convex\emph{\ }%
function on $Z$ (it might be identically $\left\{  +\infty\right\}  $).

Let $p:Z\rightarrow\mathbb{R\cup}\left\{  +\infty\right\}  $, not identically
$\left\{  +\infty\right\}  $, be given. We define its \emph{subconjugate}
$p^{\sharp}:X\rightarrow\mathbb{R\cup}\left\{  +\infty\right\}  $ by:
\begin{equation}
p^{\sharp}\left(  x\right)  =\sup_{z}\left\{  u\left(  x,z\right)  -p\left(
z\right)  \right\} \label{Fen*}%
\end{equation}

It follows from the definitions that $p^{\sharp}$ is a $u$-convex\emph{\ }%
function on $X$\ (it might be identically $\left\{  +\infty\right\}  $).

\begin{example}
Set $f\left(  x\right)  =u\left(  x,\bar{z}\right)  +a$. Then
\[
f^{\sharp}\left(  \bar{z}\right)  =\sup_{x}\left\{  u\left(  x,\bar{z}\right)
-u\left(  x,\bar{z}\right)  -a\right\}  =-a
\]

\end{example}

Conjugation reverses ordering: if $f_{1}\leq f_{2}$, then $f_{1}^{\sharp}\geq
f_{2}^{\sharp}$, and if $p_{1}\leq p_{2},$ then $p_{1}^{\sharp}\geq
p_{2}^{\sharp}$. As a consequence, if $f$ is $u$-convex, not identically
$\left\{  +\infty\right\}  $, then $f^{\sharp}$ is $u$-convex, not identically
$\left\{  +\infty\right\}  $,. Indeed, since $f$ is $u$-convex, we have
$f\left(  x\right)  \geq u\left(  x,z\right)  +a\;$for some $\left(
z,a\right)  $, and then $f^{\sharp}\left(  z\right)  \leq-a<\infty.$

\begin{proposition}
[the Fenchel inequality]For any functions $f:X\rightarrow\mathbb{R\cup
}\left\{  +\infty\right\}  $ and $p:X\rightarrow\mathbb{R\cup}\left\{
+\infty\right\}  $, not identically $\left\{  +\infty\right\}  $, we have:
\begin{align*}
\forall\left(  x,z\right)  ,\;\;f\left(  x\right)  +f^{\sharp}\left(
z\right)   & \geq u\left(  x,z\right)  \;\;\\
\forall\left(  x,z\right)  \;\;p\left(  z\right)  +p^{\sharp}\left(  x\right)
& \geq u\left(  x,z\right)  \;\;
\end{align*}

\end{proposition}

\subsection{Subgradients}

Let $f:X\rightarrow\mathbb{R\cup}\left\{  +\infty\right\}  $ be given, not
identically $\left\{  +\infty\right\}  $. Take some point $x\in X$. We shall
say that a point $z\in Z$ is a \emph{subgradient} of $f$ at $x$ if the points
$x$ and $z$ achieve equality in the Fenchel inequality:
\begin{equation}
f\left(  x\right)  +f^{\sharp}\left(  z\right)  =u\left(  x,z\right)
\label{subd}%
\end{equation}

The set of subgradients of $f$ at $x$ will be called the
\emph{subdifferential} of $f$ at $x$ and denoted by $\partial f\left(
x\right)  $.

Similarly, let $p:Z\rightarrow\mathbb{R\cup}\left\{  +\infty\right\}  $ be
given, not identically $\left\{  +\infty\right\}  $. Take some point $z\in Z$.
We shall say that a point $x\in X$ is a \emph{subgradient} of $p$ at $z$ if:
\begin{equation}
p^{\sharp}\left(  x\right)  +p\left(  z\right)  =u\left(  x,z\right)
\label{subdd}%
\end{equation}
The set of subgradients of $p$ at $z$ will be called the
$\emph{subdifferential}$ of $p$ at $z$ and denoted by $\partial p\left(
z\right)  $.

\begin{proposition}
\label{prop1}The following are equivalent:

\begin{enumerate}
\item $z\in\partial f\left(  x\right)  $

\item $\forall x^{\prime},\;\;f\left(  x^{\prime}\right)  \geq f\left(
x\right)  +u\left(  x^{\prime},z\right)  -u\left(  x,z\right)  \;\;$
\end{enumerate}

If equality holds for some $x^{\prime}$, then $z\in\partial f\left(
x^{\prime}\right)  $ as well.
\end{proposition}

\begin{proof}
We begin with proving that the first condition implies the second one. Assume
$z\in\partial f\left(  x\right)  $. Then, by (\ref{subd}) and the Fenchel
inequality, we have:
\[
f\left(  x^{\prime}\right)  \geq u\left(  x^{\prime},z\right)  -f^{\sharp
}\left(  z\right)  =u\left(  x^{\prime},z\right)  -\left[  u\left(
x,z\right)  -f\left(  x\right)  \right]
\]

We then prove that the second condition implies the first one. Using the
inequality, we have:
\begin{align*}
f^{\sharp}\left(  z\right)   & =\sup_{x^{\prime}}\left\{  u\left(  x^{\prime
},z\right)  -f\left(  x^{\prime}\right)  \right\} \\
& \leq\sup_{x^{\prime}}\left\{  u\left(  x^{\prime},z\right)  -f\left(
x\right)  -u\left(  x^{\prime},z\right)  +u\left(  x,z\right)  \right\} \\
& =u\left(  x,z\right)  -f\left(  x\right)
\end{align*}

so $f\left(  x\right)  +f^{\sharp}\left(  z\right)  \leq u\left(  x,z\right)
$. We have the converse by the Fenchel inequality, so equality holds.

Finally, if equality holds for some $x^{\prime}$ in condition (2), then
$\;f\left(  x^{\prime}\right)  -u\left(  x^{\prime},z\right)  =f\left(
x\right)  -u\left(  x,z\right)  $, so that:
\begin{align*}
\forall x^{\prime\prime},\;\;f\left(  x^{\prime\prime}\right)   & \geq
f\left(  x\right)  -u\left(  x,z\right)  +u\left(  x^{\prime\prime},z\right)
\\
& =f\left(  x^{\prime}\right)  -u\left(  x^{\prime},z\right)  +u\left(
x^{\prime\prime},z\right)
\end{align*}
which implies that $z\in\partial f\left(  x^{\prime}\right)  $.
\end{proof}

There is a similar result for functions $p:Z\rightarrow\mathbb{R\cup}\left\{
+\infty\right\}  $, not identically $\left\{  +\infty\right\}  $: we have
$x\in\partial p\left(  z\right)  $ if and only if
\begin{equation}
\forall\left(  x^{\prime},z^{\prime}\right)  ,\;\;p\left(  z^{\prime}\right)
\geq p\left(  z\right)  +u\left(  x,z^{\prime}\right)  -u\left(  x,z\right)
\;\;\label{eq5}%
\end{equation}

\subsection{Biconjugates}

It follows from the Fenchel inequality that, if $p:Z\rightarrow\mathbb{R\cup
}\left\{  +\infty\right\}  $ is not identically $\left\{  +\infty\right\}  $:%

\begin{equation}
p^{\sharp\sharp}\left(  z\right)  =\sup_{x}\left\{  u\left(  x,z\right)
-p^{\sharp}\left(  x\right)  \right\}  \leq p\left(  z\right) \label{eq7}%
\end{equation}

\begin{example}
Set $p\left(  z\right)  =u\left(  \bar{x},z\right)  +b$. Then
\begin{align*}
p^{\sharp\sharp}\left(  z\right)   & =\sup_{x}\left\{  u\left(  x,z\right)
-p^{\sharp}\left(  x\right)  \right\} \\
& \geq u\left(  \bar{x},z\right)  -p^{\sharp}\left(  \bar{x}\right) \\
& =u\left(  \bar{x},z\right)  +b=p\left(  z\right)
\end{align*}

\end{example}

This example generalizes to all $u$-convex functions. Denote by $C_{u}\left(
Z\right)  $ the set of all $u$-convex functions on $Z$.

\begin{proposition}
\label{prop2}For every function $p:Z\rightarrow\mathbb{R\cup}\left\{
+\infty\right\}  $, not identically $\left\{  +\infty\right\}  $, we have
\[
p^{\sharp\sharp}\left(  z\right)  =\sup_{\varphi}\left\{  \varphi\left(
z\right)  \;\left|  \;\varphi\leq p,\;\varphi\in\mathbb{C}_{u}\left(
Z\right)  \right.  \right\}
\]

\end{proposition}

\begin{proof}
Denote by $\bar{p}\left(  z\right)  $ the right-hand side of the above
formula. We want to show that $p^{\sharp\sharp}\left(  z\right)  =\bar
{p}\left(  z\right)  $

Since $p^{\sharp\sharp}\leq p$ and $p^{\sharp\sharp}$ is $u$-convex, we must
have$\;p^{\sharp\sharp}\leq\bar{p}$.

On the other hand, $\bar{p}$ is $u$-convex because it is a supremum of
$u$-convex functions. So there must be some $B\subset X\times\mathbb{R}$ such
that:
\[
p\left(  z\right)  =\sup_{\left(  x,b\right)  \in B}\left\{  u\left(
x,z\right)  +b\right\}
\]
Let $\left(  x,b\right)  $ $\in B$. Since $\bar{p}\leq p$, we have $u\left(
x,z\right)  +b\leq\bar{p}\left(  z\right)  \leq p\left(  z\right)  $. Taking
biconjugates, as in the preceding example, we get $\ u\left(  x,z\right)
+b\leq p^{\sharp\sharp}\left(  z\right)  $. Taking the supremums over $\left(
x,b\right)  $ $\in B$, we get the desired result.
\end{proof}

\begin{corollary}
\label{cor2}Let $p:Z\rightarrow\mathbb{R\cup}\left\{  +\infty\right\}  $ be a
$u$-convex function, not identically $\left\{  +\infty\right\}  $. Then
$p=p^{\sharp\sharp}$, and the following are equivalent:

\begin{enumerate}
\item $x\in\partial p\left(  z\right)  $

\item $p\left(  z\right)  +p^{\sharp}\left(  x\right)  =u\left(  x,z\right)  $

\item $z\in\partial p^{\sharp}\left(  x\right)  $
\end{enumerate}
\end{corollary}

\begin{proof}
We have $p^{\sharp\sharp}\leq p$ always by relation (\ref{eq7}). Since $p$ is
$u$-convex, we have:
\[
p\left(  z\right)  =\sup_{\left(  x,b\right)  \in B}\left\{  u\left(
x,z\right)  +b\right\}
\]
for some $B\subset X\times\mathbb{R}$. By proposition \ref{prop2}, we have:
\[
\sup_{\left(  x,b\right)  \in B}\left\{  u\left(  x,z\right)  +b\right\}  \leq
p^{\sharp\sharp}\left(  z\right)
\]
and so we must have $p=p^{\sharp\sharp}$. Taking this relation into account,
as well as the definition of the subgradient, we see that condition (2) is
equivalent both to (1) and to (2)
\end{proof}

\begin{definition}
We shall say that a function $p:Z\rightarrow\mathbb{R\cup}\left\{
+\infty\right\}  $ is $u$-adapted if it is not identically $\left\{
+\infty\right\}  $ and there is some $\left(  x,b\right)  \in X\times R$ such
that:
\[
\forall z\in Z,\;\;p\left(  z\right)  \geq u\left(  x,z\right)  +b
\]

\end{definition}

It follows from the above that if $p$ is $u$-adapted, then so are $p^{\sharp}%
$, $p^{\sharp\sharp}$ and all further subconjugates. Note that a $u$-convex
function which is not identically $\left\{  +\infty\right\}  $ is $u$-adapted.

\begin{corollary}
\label{cor3}Let $p:Z\rightarrow\mathbb{R\cup}\left\{  +\infty\right\}  $ be
$u$-adapted. Then :
\[
p^{\sharp\sharp\sharp}=p^{\sharp}%
\]

\end{corollary}

\begin{proof}
If $p$ is $u$-adapted, then $p^{\sharp}$ is $u$-convex and not identically
$\left\{  +\infty\right\}  $. The result then follows from corollary
\ref{cor2}.
\end{proof}

\subsection{Smoothness}

Since $u$ is continuous and $X\times Z$ is compact, the family $\left\{
u\left(  x,\cdot\right)  \;\left|  \;x\in X\right.  \right\}  $ is uniformly
equicontinuous on $Z$. It follows from the definition \ref{u-conv*} that all
$u$-convex functions on $Z$ are continuous (in particular, they are finite everywhere)..

Denote by $k$ the upper bound of $\left\|  D_{x}u\left(  x,z\right)  \right\|
$ for $\left(  x,z\right)  \in X\times Z$. Since $D_{x}u$ is continuous and
$X\times Z$ is compact, we have $k<\infty$, and the functions $x\rightarrow
u\left(  x,z\right)  $ are all $k$-Lipschitzian on $X$. Again, it follows from
the definition \ref{u-conv} that all $u$-convex functions on $X$ are
$k$-Lipschitz (in particular, they are finite everywhere). By a theorem of
Rademacher, they are differentiable almost everywhere with respect to the
Lebesgue measure..

Let $f$ $:X\rightarrow R$ be convex. Since $f=f^{\sharp\sharp}$, we have:
\[
f\left(  x\right)  =\sup_{z}\left\{  u\left(  x,z\right)  -f^{\sharp}\left(
z\right)  \right\}
\]
Since $f^{\sharp}$ is $u$-convex, it is continuous, and the supremum is
achieved on the right-hand side, at some point $z\in\partial f\left(
x\right)  $. This means that all $u$-convex functions on $X$ are
subdifferentiable everywhere on $X$.

Let $x$ be a point where $f$ is differentiable, with derivative $D_{x}f\left(
x\right)  $, and let $z\in\partial f\left(  x\right)  $. Consider the function
$\varphi\left(  x^{\prime}\right)  =u\left(  x^{\prime},z\right)  -f^{\sharp
}\left(  z\right)  $. By proposition \ref{prop1}, we have $\varphi\leq f$ and
$\varphi\left(  x\right)  =f\left(  x\right)  $, so that $\varphi$ and $f$
must have the same derivative at $x:$
\begin{equation}
D_{x}f\left(  x\right)  =D_{x}u\left(  x,z\right) \label{eq10}%
\end{equation}

By assumption (see (\ref{eq11}), this equation defines $z$ uniquely. We shall
denote it by $z=\nabla_{u}f\left(  x\right)  $. In other words, at every point
$x$ where $f$ is differentiable, the subdifferential $\partial f\left(
x\right)  $ reduces to a singleton, namely $\left\{  \nabla_{u}f\left(
x\right)  \right\}  $. Combining all this information, we get:

\begin{proposition}
\label{prop5}For every $u$-convex function $f:X\rightarrow R$, there is a map
$\nabla_{u}f:X\rightarrow Z$ such that, for almost every $x:$
\[
D_{x}f\left(  x\right)  =D_{x}u\left(  x,z\right)  \text{ }\Longleftrightarrow
z=\nabla_{u}f\left(  x\right)
\]

\end{proposition}

The following result will also be useful:

\begin{proposition}
\label{prop7} Let $p:Z\rightarrow R$ be $u$-adapted, and let $x\in X$ be
given. Then there is some point $z\in\partial p^{\sharp}\left(  x\right)  $
such that $p\left(  z\right)  =p^{\sharp\sharp}\left(  z\right)  $.
\end{proposition}

\begin{proof}
Assume otherwise, so that for every $z\in\partial p^{\sharp}\left(  x\right)
$ we have $p^{\sharp\sharp}\left(  z\right)  <p\left(  z\right)  $. For every
$z\in\partial p^{\sharp}\left(  x\right)  $, we have $x\in\partial
p^{\sharp\sharp}\left(  z\right)  $, so that, by proposition \ref{prop1}, we
have
\[
p^{\sharp\sharp}\left(  z^{\prime}\right)  \geq u\left(  x,z^{\prime}\right)
-u\left(  x,z\right)  +p^{\sharp\sharp}\left(  z\right)
\]
for all $z^{\prime}\in Z$, the inequality being strict if $z^{\prime}\notin$
$\partial p^{\sharp}\left(  x\right)  .$\ Set $\varphi_{z}\left(  z^{\prime
}\right)  =u\left(  x,z^{\prime}\right)  -u\left(  x,z\right)  +p^{\sharp
\sharp}\left(  z\right)  $. We have:
\begin{align*}
z^{\prime}  & \notin\partial p^{\sharp}\left(  x\right)  \Longrightarrow
\varphi_{z}\left(  z^{\prime}\right)  <p^{\sharp\sharp}\left(  z^{\prime
}\right)  \leq p\left(  z^{\prime}\right) \\
z^{\prime}  & \in\partial p^{\sharp}\left(  x\right)  \Longrightarrow
\varphi_{z}\left(  z^{\prime}\right)  \leq p^{\sharp\sharp}\left(  z^{\prime
}\right)  <p\left(  z^{\prime}\right)
\end{align*}
so that $\varphi_{z}\left(  z^{\prime}\right)  <p\left(  z^{\prime}\right)  $
for all $\left(  z,z^{\prime}\right)  $. Since $Z$ is compact, there is some
$\varepsilon>0$ such that $\varphi_{z}\left(  z^{\prime}\right)
+\varepsilon\leq p\left(  z^{\prime}\right)  $ for all $\left(  z,z^{\prime
}\right)  $. Taking the subconjugate with respect to $z^{\prime}$, we get:
\begin{align*}
p^{\sharp}\left(  x\right)   & \leq\sup_{z^{\prime}}\left\{  u\left(
x,z^{\prime}\right)  -\varphi_{z}\left(  z^{\prime}\right)  \right\}
-\varepsilon\\
& =\sup_{z^{\prime}}\left\{  u\left(  x,z^{\prime}\right)  -u\left(
x,z^{\prime}\right)  +u\left(  x,z\right)  -p^{\sharp\sharp}\left(  z\right)
\right\}  -\varepsilon\\
& =u\left(  x,z\right)  -p^{\sharp\sharp}\left(  z\right)  -\varepsilon
=p^{\sharp}\left(  x\right)  -\varepsilon
\end{align*}
which is a contradiction. The result follows
\end{proof}

\begin{corollary}
\label{cor7}If $x$ is a point where $p^{\sharp}$ is differentiable, then:
\begin{equation}
p\left(  \nabla_{u}p^{\sharp}\left(  x\right)  \right)  =p^{\sharp\sharp
}\left(  \nabla_{u}p^{\sharp}\left(  x\right)  \right) \label{eq25}%
\end{equation}
and:
\begin{equation}
p^{\sharp}\left(  x\right)  =u\left(  x,\nabla_{u}p^{\sharp}\left(  x\right)
\right)  -p\left(  \nabla_{u}p^{\sharp}\left(  x\right)  \right) \label{eq26}%
\end{equation}

\end{corollary}

\begin{proof}
Just apply the preceding proposition, bearing in mind that $\partial
p^{\sharp}\left(  x\right)  $ contains only one point, namely $\nabla
_{u}p^{\sharp}\left(  x\right)  $. This yields equation (\ref{eq25}) Equation
(\ref{eq26}) follows from the definition of the subgradient:
\[
p^{\sharp}\left(  x\right)  =u\left(  x,\nabla_{u}p^{\sharp}\left(  x\right)
\right)  -p^{\sharp\sharp}\left(  \nabla_{u}p^{\sharp}\left(  x\right)
\right)
\]
and equation (\ref{eq25}).
\end{proof}

\subsection{$v$-concave functions.}

Let us now consider the duality between $Y$ and $\dot{Z}$. Given $v:Y\times
Z\rightarrow R$, we say that a map $g:Y\rightarrow\mathbb{R\cup}\left\{
-\infty\right\}  $ is $v$\emph{-concave} iff there exists a non-empty subset
$A\subset Z\times\mathbb{R}$ such that:
\begin{equation}
\forall y\in Y,\;\;g\left(  y\right)  =\inf_{\left(  z,\alpha\right)  \in
A}\left\{  v\left(  y,z\right)  +a\right\}
\end{equation}
and a function $p:Z\rightarrow\mathbb{R\cup}\left\{  -\infty\right\}  $ will
be called $v$\emph{-concave} iff there exists a non-empty subset $B\subset
X\times\mathbb{R}$ such that:
\begin{equation}
p\left(  z\right)  =\inf_{\left(  x,b\right)  \in B}\left\{  v\left(
y,z\right)  +b\right\}
\end{equation}

All the results on $u$-convex functions carry over to $v$-concave functions,
with obvious modifications. The \emph{superconjugate} of a function
$g:Y\rightarrow\mathbb{R\cup}\left\{  -\infty\right\}  $, not identically
$\left\{  -\infty\right\}  $, is defined by:
\begin{equation}
g^{\flat}\left(  z\right)  =\inf_{y}\left\{  v\left(  y,z\right)  -g\left(
y\right)  \right\}
\end{equation}
and the \emph{superconjugate} of a function $p:Z\rightarrow\mathbb{R\cup
}\left\{  -\infty\right\}  $, not identically $\left\{  -\infty\right\}  $, is
given by:
\begin{equation}
p^{\flat}\left(  y\right)  =\inf_{z}\left\{  v\left(  y,z\right)  -p\left(
z\right)  \right\}
\end{equation}

The superdifferential $\partial p^{\flat}$ is defined by:
\[
\partial p^{\flat}\left(  y\right)  =\arg\min_{z}\left\{  v\left(  y,z\right)
-p\left(  z\right)  \right\}
\]
and we have the Fenchel inequality:
\[
p\left(  z\right)  +p^{\flat}\left(  y\right)  \leq v\left(  y,z\right)
\;\;\forall\left(  y,z\right)
\]
with equality iff $z\in\partial p^{\flat}\left(  y\right)  $. Note finally
that $p^{\flat\flat}\geq p$, with equality if $p$ is $v$-concave

\section{The dual optimization problem.}

Denote by $\mathcal{A}$ the set of all bounded function on $Z$:%

\[
p\in\mathcal{A}\Longleftrightarrow\sup_{z}\left|  p\left(  z\right)  \right|
<\infty
\]
and consider the minimization problem: \
\begin{equation}
\inf_{p\in\mathcal{A}}\left[  \int_{X}p^{\sharp}\left(  x\right)  d\mu
-\int_{Y}p^{\flat}\left(  y\right)  d\nu\right] \label{P}%
\end{equation}

\begin{proposition}
\label{prop6}The minimum is attained in problem (P)
\end{proposition}

Take a minimizing sequence $p_{n}$:
\[
\int_{X}p_{n}^{\sharp}\left(  x\right)  d\mu-\int_{Y}p_{n}^{\flat}\left(
y\right)  d\nu\rightarrow\inf\left(  \mathrm{P}\right)
\]
Setting $q_{n}=p_{n}+a$, for some constant $a$. Then $q_{n}^{\sharp}%
=p_{n}^{\sharp}-a$ and $q_{n}^{\flat}=p_{n}^{\flat}-a$. Since $\mu\left(
X\right)  =\nu\left(  Y\right)  $, we have:
\[
\int_{X}q_{n}^{\sharp}\left(  x\right)  d\mu-\int_{Y}q_{n}^{\flat}\left(
y\right)  d\nu=\int_{X}p_{n}^{\sharp}\left(  x\right)  d\mu-\int_{Y}%
p_{n}^{\flat}\left(  y\right)  d\nu\rightarrow\inf\left(  \mathrm{P}\right)
\]

Setting $a=-\inf_{z}p_{n}\left(  z\right)  $, we find $\inf_{z}q_{n}\left(
z\right)  =0$. So there is no loss of generality in assuming that:
\[
\forall z,\;\;\inf_{z}\;p_{n}\left(  z\right)  =0
\]
which we shall do from now on.

The sequences $p_{n}^{\sharp}$ is $k$-Lipschitzian. Since $p_{n}\geq0$, we
have
\[
p_{n}^{\sharp}\left(  x\right)  \leq\sup_{z}u\left(  x,z\right)  \leq
\max_{x,z}u\left(  x,z\right)
\]

Choose $z_{n}$ such that $p_{n}\left(  z_{n}\right)  \leq1$. We then have:
\[
p_{n}^{\sharp}\left(  x\right)  \geq u\left(  x,z_{n}\right)  -p_{n}\left(
z_{n}\right)  \geq\min_{x,z}u\left(  x,z\right)  -1
\]

So the sequence $p_{n}^{\sharp}$ is uniformly bounded. By Ascoli's theorem,
there is a uniformly convergent subsequence. Similarly, after extracting this
first subsequence, we extract another one along which $p_{n}^{\flat}$
converges uniformly. The resulting subsequence will still be denoted by
$p_{n}$, so that:
\begin{align*}
p_{n}^{\sharp}  & \rightarrow f\;\;\mathrm{uniformly}\\
p_{n}^{\flat}  & \rightarrow g\;\;\mathrm{uniformly}%
\end{align*}

Taking limits, we get:
\begin{equation}
\int_{X}f\left(  x\right)  d\mu-\int_{Y}g\left(  y\right)  d\nu=\inf\left(
\mathrm{P}\right) \label{o2}%
\end{equation}

It is easy to see that $p_{n}^{\sharp\sharp}\rightarrow f^{\sharp}$ and
$p_{n}^{\flat\flat}\rightarrow g^{\flat}$ uniformly. Since $p_{n}%
^{\sharp\sharp}\leq p_{n}$ and $p_{n}^{\flat\flat}\geq p_{n}$, we have
$p_{n}^{\sharp\sharp}\leq p_{n}^{\flat\flat}$, and hence $f^{\sharp}\leq
g^{\flat}$. Set:%

\[
\bar{p}=\frac{1}{2}\left(  f^{\sharp}+g^{\flat}\right)  .
\]

Since $f^{\sharp}$ and $g^{\flat}$ are Lipschitz continuous, so is $\bar{p} $.
Since $f^{\sharp}$ $\leq\bar{p}\leq g^{\flat}$, we must have $\bar{p}^{\sharp
}\leq f^{\sharp\sharp}\leq f$ and $\bar{p}^{\flat}\geq g^{\flat\flat}\geq g$.
Hence:
\[
\int_{X}\bar{p}^{\sharp}\left(  x\right)  d\mu-\int_{Y}\bar{p}^{\flat}\left(
y\right)  d\nu\leq\int_{X}f\left(  x\right)  d\mu-\int_{Y}g\left(  y\right)
d\nu=\inf\left(  \mathrm{P}\right)
\]

Since $f^{\sharp}$ $\leq\bar{p}\leq g^{\flat}$, both \ sides being continuous
functions, the function $\bar{p}$ must be bounded on $Z$, and the above
inequality shows that it is a minimizer. The proof is concluded.

\section{Proof of the main result.}

Let us now express the optimality condition in problem (P). Set:
\begin{align*}
\bar{s}\left(  x\right)   & =\nabla_{u}p^{\sharp}\left(  x\right) \\
\bar{t}\left(  y\right)   & =\nabla_{v}p^{\flat}\left(  y\right)
\end{align*}
where the gradient maps $\nabla_{u}p^{\sharp}$ and $\nabla_{v}p^{\flat}$ have
been defined in proposition \ref{prop5}.

\begin{proposition}
\label{prop8}$\bar{s}\left(  \mu\right)  =\bar{t}\left(  \nu\right)  $
\end{proposition}

\begin{proof}
We follow the argument in Carlier \ref{Carlier}. Take any continuous function
$\varphi:Z\rightarrow R$. Since $\bar{p}$ is a minimizer, we have, for any
integer $n$:
\begin{equation}
\int_{X}n\left[  \left(  \bar{p}+\frac{1}{n}\varphi\right)  ^{\sharp}-\bar
{p}^{\sharp}\right]  d\mu-\int_{Y}n\left[  \left(  \bar{p}+\frac{1}{n}%
\varphi\right)  ^{\flat}-\bar{p}^{\flat}\right]  d\nu\geq0\label{re}%
\end{equation}
We deal with the first integral. Set $\bar{p}+\frac{1}{n}\varphi=p_{n}$. Since
$p_{n}^{\sharp}$ is $u$-convex, it differentiable almost everywhere. Take a
negligible subset $X_{0}$ such that all the $p_{n},n\in N,$ and $\bar{p}$, are
differentiable at every $x\notin X_{0}$. If $x\notin X_{0}$, then $\nabla
_{u}p_{n}^{\sharp}\left(  x\right)  $ is the only point in $\partial
p_{n}^{\sharp}\left(  x\right)  $, and we have, by corollary \ref{cor7}:
\[
u\left(  x,\nabla_{u}p_{n}^{\sharp}\left(  x\right)  \right)  -\bar{p}\left(
\nabla_{u}p_{n}^{\sharp}\left(  x\right)  \right)  \leq\bar{p}^{\sharp}\left(
x\right)  =u\left(  x,\nabla_{u}\bar{p}^{\sharp}\left(  x\right)  \right)
-\bar{p}\left(  \nabla_{u}\bar{p}^{\sharp}\left(  x\right)  \right)
\]
so that:
\begin{equation}
u\left(  x,\nabla_{u}p_{n}^{\sharp}\left(  x\right)  \right)  -\bar{p}\left(
\nabla_{u}p_{n}^{\sharp}\left(  x\right)  \right)  -u\left(  x,\nabla_{u}%
\bar{p}^{\sharp}\left(  x\right)  \right)  +\bar{p}\left(  \nabla_{u}\bar
{p}^{\sharp}\left(  x\right)  \right)  \leq0\label{ql1}%
\end{equation}
yielding:
\begin{equation}
p_{n}^{\sharp}\left(  x\right)  +\frac{1}{n}\varphi\left(  \nabla_{u}%
p_{n}^{\sharp}\left(  x\right)  \right)  -\bar{p}^{\sharp}\left(  x\right)
\leq0\label{eq50}%
\end{equation}
From the definition of $p_{n}^{\sharp}$, we have, using corollary \ref{cor7}
again:
\[
u\left(  x,\nabla_{u}\bar{p}^{\sharp}\left(  x\right)  \right)  -\bar
{p}\left(  \nabla_{u}\bar{p}^{\sharp}\left(  x\right)  \right)  -\frac{1}%
{n}\varphi\left(  \nabla_{u}\bar{p}^{\sharp}\left(  x\right)  \right)  \leq
p_{n}^{\sharp}\left(  x\right)  =u\left(  x,\nabla_{u}p_{n}^{\sharp}\left(
x\right)  \right)  -\bar{p}\left(  \nabla_{u}p_{n}^{\sharp}\left(  x\right)
\right)  -\frac{1}{n}\varphi\left(  \nabla_{u}p_{n}^{\sharp}\left(  x\right)
\right)
\]
Rewriting this, we get:
\begin{align}
\frac{1}{n}\varphi\left(  \nabla_{u}p_{n}^{\sharp}\left(  x\right)  \right)
-\frac{1}{n}\varphi\left(  \nabla_{u}\bar{p}^{\sharp}\left(  x\right)
\right)   & \leq u\left(  x,\nabla_{u}p_{n}^{\sharp}\left(  x\right)  \right)
-\bar{p}\left(  \nabla_{u}p_{n}^{\sharp}\left(  x\right)  \right)  -u\left(
x,\nabla_{u}\bar{p}^{\sharp}\left(  x\right)  \right)  +\bar{p}\left(
\nabla_{u}\bar{p}^{\sharp}\left(  x\right)  \right) \label{ql2}\\
& =p_{n}^{\sharp}\left(  x\right)  +\frac{1}{n}\varphi\left(  \nabla_{u}%
p_{n}^{\sharp}\left(  x\right)  \right)  -\bar{p}^{\sharp}\left(  x\right)
\label{ql3}%
\end{align}

yielding:
\begin{equation}
-\frac{1}{n}\varphi\left(  \nabla_{u}\bar{p}^{\sharp}\left(  x\right)
\right)  \leq p_{n}^{\sharp}\left(  x\right)  -\bar{p}^{\sharp}\left(
x\right) \label{eq51}%
\end{equation}
Now let $n\rightarrow\infty$. Using corollary \ref{cor7}, we have:
\[
u\left(  x,\nabla_{u}p_{n}^{\sharp}\left(  x\right)  \right)  =p_{n}\left(
\nabla_{u}p_{n}^{\sharp}\left(  x\right)  \right)  +p_{n}^{\sharp}\left(
x\right)
\]
Since $Z$ is compact, the sequence $\nabla_{u}p_{n}^{\sharp}\left(  x\right)
\in Z$ has a cluster point $z,$ and since $p_{n}$ and $p_{n}^{\sharp}$
converge to $\bar{p}$ and $\bar{p}^{\sharp}$ uniformly, we get in the limit:
\[
u\left(  x,z\right)  =\bar{p}\left(  z\right)  +\bar{p}^{\sharp}\left(
x\right)
\]
so that $z\in\partial\bar{p}^{\sharp}\left(  x\right)  $. But that
subdifferential consists only of the point $\nabla_{u}\bar{p}^{\sharp}\left(
x\right)  $, so that $z=\nabla_{u}\bar{p}^{\sharp}\left(  x\right)  $. This
shows that the cluster point $z$ is unique, so that the whole sequence must
converge:
\[
\nabla_{u}p_{n}^{\sharp}\left(  x\right)  \rightarrow\nabla_{u}\bar{p}%
^{\sharp}\left(  x\right)
\]

Inequalities () and () together give:
\begin{equation}
-\varphi\left(  \nabla_{u}\bar{p}^{\sharp}\left(  x\right)  \right)  \leq
n\left(  p_{n}^{\sharp}\left(  x\right)  -\bar{p}^{\sharp}\left(  x\right)
\right)  \leq\varphi\left(  \nabla_{u}p_{n}^{\sharp}\left(  x\right)  \right)
\label{eq52}%
\end{equation}
Taking limits in the inequalities (\ref{eq50}) and (\ref{eq51}), we get:
\[
\forall x\notin X_{0},\;\;\lim_{n}\;n\left(  p_{n}^{\sharp}\left(  x\right)
-\bar{p}^{\sharp}\left(  x\right)  \right)  =-\varphi\left(  \nabla_{u}\bar
{p}^{\sharp}\left(  x\right)  \right)
\]

Similarly, we have:
\[
\forall y\notin Y_{0},\;\;\lim_{n}\;n\left(  p_{n}^{\flat}\left(  x\right)
-\bar{p}^{\flat}\left(  x\right)  \right)  =-\varphi\left(  \nabla_{u}\bar
{p}^{\flat}\left(  x\right)  \right)
\]
where $Y_{0}\subset Y$ is negligible.

Because of (\ref{eq52}), we can apply the dominated convergence theorem to
inequality (\ref{re}). We get:
\[
-\int_{X}\varphi\left(  \nabla\bar{p}^{\sharp}\left(  x\right)  \right)
d\mu+\int_{Y}\varphi\left(  \nabla\bar{p}^{\flat}\left(  y\right)  \right)
d\nu\geq0
\]
Since the inequality must hold for $-\varphi$ as well as $\varphi$, it is in
fact an equality. In other words, for any $\varphi:Z\rightarrow R$ with
compact support, we have:
\begin{equation}
\int_{X}\varphi\left(  \nabla\bar{p}^{\sharp}\left(  x\right)  \right)
d\mu=\int_{Y}\varphi\left(  \nabla\bar{p}^{\flat}\left(  y\right)  \right)
d\nu\label{ef1}%
\end{equation}
and this means that $\bar{s}\left(  \mu\right)  =\bar{t}\left(  \nu\right)  $,
as announced.
\end{proof}

Set $\bar{s}\left(  \mu\right)  =\bar{t}\left(  \nu\right)  =\lambda$. This is
a positive measure on $Z$, not necessarily absolutely continous with respect
to the Lebesgue measure.

Applying corollary \ref{cor7}, we have:
\begin{align*}
\bar{p}^{\sharp}\left(  x\right)   & =u\left(  x,\bar{s}\left(  x\right)
\right)  -\bar{p}\left(  \bar{s}\left(  x\right)  \right) \\
\bar{p}^{\flat}\left(  y\right)   & =v\left(  y,\bar{t}\left(  y\right)
\right)  -\bar{p}\left(  \bar{t}\left(  y\right)  \right)
\end{align*}
and hence:
\begin{align}
\int_{X}\bar{p}^{\sharp}\left(  x\right)  d\mu-\int_{Y}\bar{p}^{\flat}\left(
y\right)  d\nu & =\int_{X}\left[  u\left(  x,\bar{s}\left(  x\right)  \right)
-\bar{p}\left(  \bar{s}\left(  x\right)  \right)  \right]  d\mu-\int
_{Y}\left[  v\left(  y,\bar{t}\left(  y\right)  \right)  -\bar{p}\left(
\bar{t}\left(  y\right)  \right)  \right] \nonumber\\
& =\int_{X}u\left(  x,\bar{s}\left(  x\right)  \right)  d\mu-\int_{Y}v\left(
y,\bar{t}\left(  y\right)  \right)  d\nu-\left[  \int_{X}\bar{p}\left(
\bar{s}\left(  x\right)  \right)  d\mu-\int_{Y}\bar{p}\left(  \bar{t}\left(
y\right)  \right)  d\nu\right] \nonumber\\
& =\int_{X}u\left(  x,\bar{s}\left(  x\right)  \right)  d\mu-\int_{Y}v\left(
y,\bar{t}\left(  y\right)  \right)  d\nu-\left[  \int_{Z}\bar{p}\left(
z\right)  d\lambda-\int_{Z}\bar{p}\left(  z\right)  d\lambda\right]
\nonumber\\
& =\int_{X}u\left(  x,\bar{s}\left(  x\right)  \right)  d\mu-\int_{Y}v\left(
y,\bar{t}\left(  y\right)  \right)  d\nu\label{eq30}%
\end{align}
Let $\left(  s,t\right)  $ be a pair of Borelian maps such that $s\left(
\mu\right)  =t\left(  \nu\right)  $. Then, by the Fenchel inequality::%

\begin{align*}
\int_{X}u\left(  x,s\left(  x\right)  \right)  d\mu-\int_{Y}v\left(
y,t\left(  y\right)  \right)  d\nu & \leq\int_{X}\left[  \bar{p}^{\sharp
}\left(  x\right)  +\bar{p}\left(  s\left(  x\right)  \right)  \right]
d\mu-\int_{Y}\left[  p^{\flat}\left(  y\right)  +\bar{p}\left(  t\left(
y\right)  \right)  \right]  d\nu\\
& =\int_{X}\bar{p}^{\sharp}\left(  x\right)  d\mu-\int_{Y}\bar{p}^{\flat
}\left(  y\right)  d\nu+\left[  \int_{X}\bar{p}\left(  s\left(  x\right)
\right)  d\mu-\int_{Y}\bar{p}\left(  t\left(  y\right)  \right)  d\nu\right]
\end{align*}

The last bracket vanishes because $s\left(  \mu\right)  =t\left(  \nu\right)
$. Applying inequality \ref{eq30}, we get:
\[
\int_{X}u\left(  x,s\left(  x\right)  \right)  d\mu-\int_{Y}v\left(
y,t\left(  y\right)  \right)  d\nu\leq\int_{X}u\left(  x,\bar{s}\left(
x\right)  \right)  d\mu-\int_{Y}v\left(  y,\bar{t}\left(  y\right)  \right)
d\nu
\]

This shows that $\left(  \bar{s},\bar{t}\right)  $ is a maximizer,and proves
the existence part of theorem \ref{thm1}.

As for uniqueness, assume that there is another maximizer $\left(  s^{\prime
},t^{\prime}\right)  $. The preceding inequality then becomes an equality:
\[
\int_{X}u\left(  x,s\left(  x\right)  \right)  d\mu-\int_{Y}v\left(
y,t\left(  y\right)  \right)  d\nu=\int_{X}\left[  \bar{p}^{\sharp}\left(
x\right)  +\bar{p}\left(  s\left(  x\right)  \right)  \right]  d\mu-\int
_{Y}\left[  p^{\flat}\left(  y\right)  +\bar{p}\left(  t\left(  y\right)
\right)  \right]  d\nu
\]
which we rewrite as:
\[
\int_{X}\left[  \bar{p}^{\sharp}\left(  x\right)  +\bar{p}\left(  s\left(
x\right)  \right)  -u\left(  x,s\left(  x\right)  \right)  \right]  d\mu
+\int_{Y}\left[  v\left(  y,t\left(  y\right)  \right)  -p^{\flat}\left(
y\right)  +\bar{p}\left(  t\left(  y\right)  \right)  \right]  d\nu=0
\]

Both integrands are non-negative by the Fenchel inequality. If the sum is
zero, each integral must vanish, and since the integrands are non-negative,
each integrand must vanish almost everywhere. This means that:%

\begin{align*}
s\left(  x\right)   & \in\partial\bar{p}^{\sharp}\left(  x\right)
\;\;\mathrm{a.e.}\\
t\left(  y\right)   & \in\partial\bar{p}^{\flat}\left(  y\right)
\;\;\mathrm{a.e.}%
\end{align*}
and since $\partial\bar{p}^{\sharp}\left(  x\right)  =\left\{  \bar{s}\left(
x\right)  \right\}  $ and $\partial\bar{p}^{\flat}\left(  y\right)  =\left\{
\bar{t}\left(  y\right)  \right\}  $ almost everywhere, the result is proved.

\section{Some consequences.}

We shall now investigate some properties of the function $\bar{p}:Z\rightarrow
R$.

Recall that we denote $\lambda=\bar{s}\left(  \mu\right)  =\bar{t}\left(
\nu\right)  $. It is a positive measure on $Z$. Its support $\mathrm{Supp}%
\left(  \lambda\right)  $ is the complement of the largest open subset
$\Omega\subset Z$ such that $\lambda=0$ on $\Omega$. If for instance $\mu$ and
$\nu$ have the property that the measure of any open non-empty subset is
positive, then:
\[
\mathrm{Supp}\left(  \lambda\right)  =\bar{s}\left(  X\right)  =\bar{t}\left(
Y\right)
\]

\begin{proposition}
\label{prop9}$\bar{p}^{\sharp\sharp}=\bar{p}=\bar{p}^{\flat\flat}$ on
$\mathrm{Supp}\left(  \lambda\right)  $
\end{proposition}

\begin{proof}
We have seen that $\bar{p}^{\sharp\sharp}\left(  \bar{s}\left(  x\right)
\right)  =\bar{p}\left(  \bar{s}\left(  x\right)  \right)  $ $\mu$-almost
everywhere. Since $\bar{s}\left(  \mu\right)  =\lambda$, this means that
$\bar{p}^{\sharp\sharp}\left(  z\right)  =\bar{p}\left(  z\right)  $ $\lambda
$-almost everywhere. Since $\bar{p}$ and $\bar{p}^{\sharp\sharp}$ are
continuous, equality extends to the support of $\lambda$. Similarly, $\bar
{p}=\bar{p}^{\flat\flat}$ on $\mathrm{Supp}\left(  \lambda\right)  $, and the
result follows
\end{proof}

Let us illustrate this with an example

\subsection{The linear case.}

Suppose $X,Y,Z$ are compact subsets of $R^{n}$, and we want to minimize:
\begin{equation}
\int_{X}\left\|  x-s\left(  x\right)  \right\|  ^{2}d\mu-\int_{Y}\left\|
y-t\left(  y\right)  \right\|  ^{2}d\nu\label{55}%
\end{equation}
among all maps $\left(  s,t\right)  $ such that $s\left(  \mu\right)
=t\left(  \nu\right)  $. Developing the squares, this amounts to minimizing:
\[
\left[  \int_{X}x^{2}d\mu-\int_{Y}y^{2}d\nu\right]  +\left[  \int_{X}s\left(
x\right)  ^{2}d\mu-\int_{Y}t\left(  y\right)  ^{2}d\nu\right]  -\left[
\int_{X}x^{\prime}s\left(  x\right)  d\mu-\int_{Y}y^{\prime}t\left(  y\right)
d\nu\right]
\]

The first bracket is a constant (it does not depend on the choice of $s$ and
$t$). The second bracket vanishes because $s\left(  \mu\right)  =t\left(
\nu\right)  $. We are left wiht the last one. So the problem amounts to
maximizing:
\[
\int_{X}x^{\prime}s\left(  x\right)  d\mu-\int_{Y}y^{\prime}t\left(  y\right)
d\nu
\]
and it falls within the scope of theorem \ref{thm1} by setting $u\left(
x,z\right)  =x^{\prime}z$ and $v\left(  y,z\right)  =y^{\prime}z$. Then
$u$-convex functions are convex in the usual sense, $v$-concave functions are
concave in the usual sense. By proposition \ref{prop9}, $p$ is linear on
$\mathrm{Supp}\left(  \lambda\right)  $, so we may take:
\[
p\left(  z\right)  =\pi^{\prime}z
\]
for some vector $\pi\in R^{n}$. We then get $\bar{s}\left(  x\right)  $ by
maximizing $\left(  x-\pi\right)  ^{\prime}z$ over $Z$. Similarly, we get
$\bar{t}\left(  y\right)  $ by minimizing $\left(  y-\pi\right)  ^{\prime}z$
over $Z$. Note that this implies that $\bar{s}\left(  X\right)  =\bar
{t}\left(  Y\right)  \subset\partial Z$, the boundary of $Z$. Let us summarize:

\begin{proposition}
There is a single map $\left(  \bar{s},\bar{t}\right)  $ which minimizes the
integral \ref{eq55} among all maps $\left(  s,t\right)  $ such that $s\left(
\mu\right)  =t\left(  \nu\right)  $. It is given by:
\begin{align}
\bar{s}\left(  x\right)   & =\arg\max_{z}\left(  x-\pi\right)  ^{\prime
}z\label{56}\\
\bar{t}\left(  y\right)   & =\arg\min_{z}\left(  y-\pi\right)  ^{\prime
}z\label{57}%
\end{align}
and the actual value of $\pi\in R^{n}$ is found by substituting (\ref{56}) and
(\ref{57}) in the integral (\ref{55}$\ $\ and by minimizing the resulting
function of $a$.
\end{proposition}

This first example is degenerate:\ the dimension \ of $\mathrm{Supp}\left(
\lambda\right)  $ is strictly less than the dimension $n$ of $X,Y$ and $Z$.
Let us now go in the opposite direction.

\subsection{The non-degenerate case.}

Suppose $X,Y,Z$ are compact subsets of $R^{n}$. Let $\bar{z}$ belong to the
interior of $\mathrm{Supp}\left(  \lambda\right)  $, so that $\bar{z}=\bar
{s}\left(  \bar{x}\right)  =\bar{t}\left(  \bar{y}\right)  $and suppose there
are neignbourhoods $\Omega_{x}$, $\Omega_{y}$, and $\Omega_{z}$ of $\bar
{x},\bar{y}$ and $\bar{z}$ such that the restrictions $\bar{s}:\Omega
_{x}\rightarrow\Omega_{z}$ and $\bar{t}:\Omega_{y}\rightarrow\Omega_{z}$ are
invertible, with continuous inverses $\sigma:\Omega_{z}\rightarrow\Omega_{x}$
and $\ \tau:\Omega_{z}\rightarrow\Omega_{y}$.

The function $\bar{p}$ then satisfies two partial differential equations on
$\Omega_{z}$: a second-order equation of Monge-Amp\`{e}re type, and a
fourth-order equation of Euler-Lagrange type

\subsubsection{A Monge-Amp\`{e}re equation.}

.Write the definitions of $\sigma\left(  z\right)  $ and $\tau\left(
z\right)  $, for $z\in\Omega_{z}$:
\begin{align}
D_{z}u\left(  \sigma\left(  z\right)  ,z\right)   & =D_{z}\bar{p}\left(
z\right) \label{eq71}\\
D_{z}v\left(  \tau\left(  z\right)  ,z\right)   & =D_{z}\bar{p}\left(
z\right) \label{eq72}%
\end{align}

Inverting the first equation expresses $\sigma\left(  z\right)  $ in terms of
$D_{z}p\left(  z\right)  $. Inverting the second one expresses $\tau\left(
z\right)  $ \ in terms of $D_{z}p\left(  z\right)  $. Substituting into the
equation $\sigma^{-1}\circ\tau\left(  \mu\right)  =\nu$ gives a second-order
partial differential equation for $p$.

Let us give an example. Consider the problem of minimizing the integral:
\begin{equation}
\int_{X}\frac{\alpha}{2}\left\|  x-s\left(  x\right)  \right\|  ^{2}%
dx+\int_{Y}\frac{1}{2}\left\|  y-t\left(  y\right)  \right\|  ^{2}dy\label{60}%
\end{equation}
over all maps $\left(  s,t\right)  $ such that $s\left(  \mu\right)  =t\left(
\nu\right)  .$ Here $\alpha>0$ is a given constant. We apply theorem
\ref{thm1} with $u\left(  x,z\right)  =-\frac{\alpha}{2}\left\|  x-z\right\|
^{2}$ and $v\left(  y,z\right)  =\frac{1}{2}\left\|  y-z\right\|  ^{2}$.

Assume\ $\Omega_{x},\Omega_{y}$ and $\Omega_{z}$ are as above. Equations
(\ref{eq71}) and (\ref{eq72}) become:%

\begin{align*}
D_{z}\bar{p}\left(  z\right)   & =-\alpha\left(  \sigma\left(  z\right)
-z\right) \\
D_{z}\bar{p}\left(  z\right)   & =\left(  \tau\left(  z\right)  -z\right)
\end{align*}
So $\sigma\left(  z\right)  =z-\frac{1}{\alpha}D_{z}\bar{p}\left(  z\right)  $
and $\tau\left(  z\right)  =z+D_{z}\bar{p}\left(  z\right)  $. The map
$\sigma^{-1}\circ\tau$ sends $\mu$ on $\nu$. We must have:
\[
\left[  \det D_{z}\sigma\left(  z\right)  \right]  ^{-1}\det D_{z}\tau\left(
z\right)  =1
\]
and this gives a second-order equation for $\bar{p}:$%
\begin{equation}
\det\left[  I+D_{zz}^{2}\bar{p}\left(  z\right)  \right]  =\det\left[
I-\frac{1}{\alpha}D_{zz}^{2}\bar{p}\left(  z\right)  \right] \label{64}%
\end{equation}

\subsubsection{An Euler-Lagrange equation.}

Assume that the constraints on $\bar{p}$ are not binding on $\Omega_{z}$. In
other words, there is some $\varepsilon>0$ such that, for every $h$ such that
$\left|  h\right|  <\varepsilon$, and every $C^{\infty}$ function $\varphi$
with compact support in $\Omega_{z}$, the function $p_{h}=p+h\varphi$ is still
$u$-convex and $v$-concave.

Recall that $\bar{p}$ solves the optimization problem:%

\[
\inf_{p}\left[  \int_{X}p^{\sharp}\left(  x\right)  d\mu-\int_{Y}p^{\flat
}\left(  y\right)  d\nu\right]
\]
and this implies that:
\[
\int_{\Omega_{x}}p_{h}^{\sharp}\left(  x\right)  d\mu-\int_{\Omega_{y}}%
p_{h}^{\flat}\left(  y\right)  d\nu\geq\int_{\Omega_{x}}\bar{p}^{\sharp
}\left(  x\right)  d\mu-\int_{\Omega_{y}}\bar{p}^{\flat}\left(  y\right)  d\nu
\]
for every $h$. Expressing the sub- and superconjugates in terms of the sub-
and superdifferentials, and taking advantage of the fact that $p_{h}%
^{\sharp\sharp}=p=p_{h}^{\flat\flat}$, we get, for every $h:$%
\begin{align}
& \int_{\Omega_{z}}\left[  u\left(  \sigma_{h}\left(  z\right)  ,z\right)
-p_{h}\left(  z\right)  \right]  d\left[  \sigma_{h}^{-1}\left(  \mu\right)
\right]  -\int_{\Omega_{z}}\left[  v\left(  \tau_{h}\left(  z\right)
,z\right)  -p_{h}\left(  z\right)  \right]  d\left[  \tau_{h}^{-1}\left(
\nu\right)  \right] \\
& \geq\int_{\Omega_{z}}\left[  u\left(  \sigma\left(  z\right)  ,z\right)
-\bar{p}\left(  z\right)  \right]  d\left[  \sigma^{-1}\left(  \mu\right)
\right]  -\int_{\Omega_{z}}\left[  v\left(  \tau\left(  z\right)  ,z\right)
-\bar{p}\left(  z\right)  \right]  d\left[  \tau^{-1}\left(  \nu\right)
\right]
\end{align}
where $\sigma_{h}$ and $\tau_{h}$ are defined by:
\begin{align}
D_{z}u\left(  \sigma_{h}\left(  z\right)  ,z\right)   & =D_{z}p_{h}\left(
z\right) \\
D_{z}v\left(  \tau_{h}\left(  z\right)  ,z\right)   & =D_{z}p_{h}\left(
z\right)
\end{align}

Note that for $h=0$, we have $\sigma_{h}=\sigma$, $\tau_{h}=\tau$ and
$p_{h}=\bar{p}$, while $\sigma\left(  \mu\right)  =\tau\left(  \mu\right)  $.
Letting $h\rightarrow0$, we find that $\bar{p}$ must satisfy an Euler-Lagrange
equation of the fourth order.

Let us illustrate this with example (\ref{60}).

A function $p:Z\rightarrow R$\ is $u$-convex iff $D_{zz}^{2}p\left(  z\right)
\geq-\alpha I$ for every $z$, where $D_{zz}^{2}p\left(  z\right)  $ is the
Hessian matrix at $z$, and it is $v$-concave iff $D_{zz}^{2}p\left(  z\right)
\leq I$. So there are many functions $p:Z\rightarrow R$ which are both
$u$-convex and $v$-concave: they must satisfy $-\alpha I\leq D_{zz}%
^{2}p\left(  z\right)  \leq I$ everywhere.

Assume\ $\Omega_{x},\Omega_{y}$ and $\Omega_{z}$ are as above, and $-\alpha
I<D_{zz}^{2}\bar{p}\left(  z\right)  <I$ on $\Omega_{z}$. Then the integral to
be maximized with respect to $p=p_{h}$ is:
\begin{align*}
& \int_{\Omega_{z}}\left[  u\left(  z-\frac{1}{\alpha}D_{z}p\left(  z\right)
,z\right)  -p\left(  z\right)  \right]  \left[  I-\frac{1}{\alpha} D_{zz}%
^{2}p\left(  z\right)  \right]  ^{-1}dz\\
& -\int_{\Omega_{z}}\left[  v\left(  z+D_{z}p\left(  z\right)  ,z\right)
-p\left(  z\right)  \right]  \left[  I-D_{zz}^{2}p\left(  z\right)  \right]
^{-1}dz=\\
& -\int_{\Omega_{z}}\left(  \left[  \frac{1}{2\alpha}\left\|  D_{z}p\right\|
^{2}-p\right]  \left[  I-\frac{1}{\alpha}D_{zz}^{2}p\right]  ^{-1}+\left[
\frac{1}{2}\left\|  D_{z}p\right\|  ^{2}-p\right]  \left[  I-D_{zz}%
^{2}p\right]  ^{-1}\right)  dz
\end{align*}

So $\bar{p}$, which minimizes the last integral, must satisfy the
corresponding fourth order Euler-Lagrange equation We will not write it down
explicitly, although relation (\ref{64}), which is valid for $p=\bar{p}$,
would introduce some simplifications.

\end{document}